\documentclass{amsart}

\usepackage{fancyhdr}
\usepackage{amssymb}
\usepackage{amsmath}
\usepackage{latexsym}
\input xy 
\xyoption{all}
\newtheorem{theorem}{Theorem}
\newtheorem{teo}{Theorem}[section]
\newtheorem{fact}[teo]{Fact}
\newtheorem{prop}[teo]{Proposition}
\newtheorem{lem}[teo]{Lemma}
\newtheorem{cor}[teo]{Corollary}
\newtheorem{notadefi}[teo]{Notation and Definition}

\newtheorem{definition}[teo]{Definition}

\newtheorem{rem}[teo]{Remark}

\def\dnfo{\,\raise.2em\hbox{$\,\mathrel|\kern-.9em\lower.35em\hbox{$\smile$}$}}
\def\dfo{\;\raise.2em\hbox{$\mathrel|\kern-.9em\lower.35em\hbox{$\smile$}\kern-.7em\hbox{\char'57}$}\;}

\newcommand{\td}{tr.dg}
\newcommand{\s}{SU}

\newcommand{\si}{\sigma}
\newcommand{\na}{\mathbb{N}}

\newcommand{\gr}{\mathbb G}
\newcommand{\sse}{\subseteq}
\DeclareMathOperator{\Fix}{Fix}
\begin{document}
\title{ Zilber's Dichotomy for Differentially Closed Fields with an Automorphism}
\author{Ronald F. Bustamante Medina}
\email{ronald.bustamante@ucr.ac.cr}
\address{Escuela de Matem\'aticas,
Universidad de Costa Rica,
Sede Rodrigo Facio,
2060, San Jos\'e, Costa Rica}
\thanks{Research supported by Centro de Investigaciones en Matem\'atica Pura y Aplicada and Escuela de Matem\'atica, Universidad de Costa Rica}

\date{March, 2020}

\begin{abstract}
The theory of difference-differential fields of characteristic zero has a model-companion denoted by $\it DCFA$. In \cite{rbdcfa1} we prove its main properties. In \cite{rbjets} we proved a weak version of Zilber's dichotomy for {\it DCFA}.
In this paper we use arc spaces techniques as developed by Moosa, Pillay and Scanlon in \cite{arcs} to suppress the extra hypothesis needed in \cite{rbjets} and prove the full Zilber's dichotomy for {\it DCFA}, we also state how these techniques generalise to partial differential fields with an automorphism.
\end{abstract}

\maketitle

\section{Introduction and preliminaries}
A difference-differential field is a differential field with an automorphism which commutes with the derivation.
For the case of characteristic zero, the theory of difference-differential has a model-companion that we denote {\it DCFA}. 

If $(K, \si, D)$ is a model of {\it DCFA}, the fixed field of $K$ is $\Fix \si=\{x \in K: \si(x)=x \}$ and the
field of constants of $K$ is ${\mathcal C}=\{x \in K: Dx=0\}$.

Our goal is to prove the following theorem, known as Zilber dichotomy.

\begin{theorem}
\label{teo1}

If $p$ is a type of $\s$-rank 1, then it is either 1-based or non orthogonal to $\Fix \si \cap {\mathcal C}$.
\end{theorem}

This dichotomy states, intuitively speaking, that a "minimal" set is either "algebraic" or "geometrically simple".
Analogues for this dichotomy are satisfied in both differentially closed fields (models of {\it DCF}, the model companion
or the theory of differential fields) and
 existentially closed difference fields (models of {\it ACFA}, the model companion of the theory of difference fields).
Both proofs (\cite{hrusok}, and \cite{salinas}) depend heavily on Zariski geometries. These results lead to Hrushovsky's proofs of
the Mordell-Lang conjecture and the Manin-Mumfurd conjecture, respectively.  
Pillay and Ziegler
proved both dichotomies for the case of characteristic zero \cite{jets} using
algebraic jet spaces. A key point in their
proofs is the fact that in models of {\it DCF}   $tp(a/F)$  has finite rank if and only
if the transcendence degree of the differential field generated by the tuple $a$ and $F$ over $F$ is finite, 
where $F$ is a differential subfield of a the model.  
This transcendence degree is what we call the dimension of $a$ over $F$, denoted $\dim (a/F)$. An analogous statement
holds for models of {\it ACFA}.These  equivalences allow
the authors to replace an element of finite rank with something interdefinable which, by finite-dimensionality, 
turns out to be a finite sequence.
Unfortunately this equivalence does not hold for {\it DCFA}. In \cite{rbrank1} we give an example of a type
with $\s$-rank 1, but with infinite dimension. So following the jet spaces techniques we were able to prove the following weaker dichotomy.

\begin{theorem}
\emph{\cite{rbjets}.}
\label{J313}
Let $(\mathcal{U},\si, D)$ be a saturated model of {\it DCFA} and let $K=acl(K)\subset \mathcal{U}$. Let $a \in {\mathcal U}$ such that $\s(tp(a/K))=1$ and $tp(a/K)$ is finite-dimensional. 
Then $tp(a/K)$ is either 1-based or non orthogonal to $\Fix \si \cap {\mathcal C}$.
\end{theorem}

In this paper we want to suppress the finite-dimensional hypothesis. For this we follow \cite{arcs} where
 the authors use algebraic arc spaces to prove a dichotomy for differentially closed fields
with finitely many commuting derivations ($DCF_n$). The main theorem states  that if  $p$ is a regular non-locally modular type
then there is a definable subgroup of the additive group whose generic type is regular and non orthogonal to $p$.

We want to prove an analogue theorem for our case, namely:

\begin{theorem}\label{teoarcs}
 Let $p$ be a regular non-locally modular type in $\mathcal U$. Then there is a definable subgroup of the additive group
 whose generic type is regular and non-orthogonal to $p$.
\end{theorem}

How can we link Theorem \ref{J313} and Theorem \ref{teoarcs} to prove Theorem \ref{teo1}? The next proposition on definable subgroups
of the additive group will do the job:
\begin{prop}
\emph{\cite{rbrank1}.}
\label{arc18}
Let $G$ be a definable subgroup of $\gr_a^n$. 

\begin{enumerate}
\item  $G$ has no proper subgroup of
finite index. 
\item $G$ is quantifier-free definable.

\item If $H$ a definable subgroup of $G$. Then $G/H$ is definably isomorphic
to a subgroup of $\gr_a^l$ for some $l$.

\item If $G$ has infinite dimension then $\s(G)\geq \omega$.

\end{enumerate}
\end{prop}

Now, provided that Theorem \ref{teoarcs} holds, we can easily prove Theorem \ref{teo1}.\\

{\it Proof of  Theorem \ref{teo1}}:\\

Suppose that $p$ is not 1-based, then it is non-locally modular. As $\s(p)=1$ it is regular.
By Theorem \ref{teoarcs} there is a definable subgroup $G$ of $\gr_a$ whose generic type $q$ is regular
 and non-orthogonal to $p$.
 
 $p \not \perp q$ implies that $\s(q)=\alpha+1$ for some $\alpha$.
Then, by 5.4.3 of \cite{wag}, $G$ contains a definable subgroup $N$ 
such that $\s(G/N)<\omega$, and by \ref{arc18} $G$ must be 
finite-dimensional.
 
By Theorem \ref{J313}, $q \not\perp Fix \si \cap {\mathcal C}$
and by transitivity $p \not\perp \Fix \si \cap {\mathcal C}$.\\
$\Box$

$DCF$ is complete, $\omega$-stable and eliminates quantifiers, these facts are key in \cite{arcs}, on the other hand {\it DCFA} is not complete, no completion is stable nor eliminates quantifiers. However its completions are easily described,  they are supersimple, quantifier-free stable and eliminates imaginaries. So we need to state some facts that 
allow us to build a set-up similar to the one exposed in \cite{arcs}.
Section~\ref{sec:supersimple} is devoted to this. In section~\ref{sec:dcfa} we give a brief description and list the main theorems on {\it DCFA}.
In section~\ref{sec:arcs} we use arc spaces to give a proof of  \ref{teoarcs}. Finally, in section~\ref{sec:dncfa} we give details on how the same reasonings apply to partial differential fields with an automorphism.

\section{Types in supersimple theories}\label{sec:supersimple}

Throughout this section $T$ shall denote  a supersimple theory which eliminates imaginaries and  $M$  a saturated model of $T$.  This implies
that there is a good notion of independence, and thus of forking. Moreover all types are ranked by the $\s$-rank.
With independence we define orthogonality: if $A \subset M$  is a set of parameters and
$p, q$ are complete types over $A$. We say that $p$ and $q$ are orthogonal ($p \perp q$) if for every $B \supseteq A$ and realisations $a$ of $p$ and $b$ of $q$, we have that $a$ is independent of $b$ over $B$ 
($a \dnfo_B b$).

\begin{prop}
\emph{(\cite{wag}, 5.1.12)}
\label{wgeq}
Let  $a \in M  $ and $A \subset M  $.
Let us suppose that $SU(a/A)= \beta + \omega^{\alpha} \cdot n $, with $n >0 $ and
  $\omega^{\alpha+1} \leq \beta < \infty $ or $\beta =0  $.
Then $tp(a/A)  $ is non-orthogonal to a type of $\s$-rank $\omega^{\alpha}$. Moreover there is $b\in acl(Aa)$
with $\s(b/A)=\omega^{\alpha}n$.
\end{prop}

\begin{definition}\label{prt3}
Let $p, \, q \in S(A)$. We say that $q$ is $p$-internal if for every realization
$a$ of $q$ there is a set $B$ such that $B \dnfo_A a$ and a tuple $c$
of realizations of $p$ such that $a \in dcl(Bc)$.
A set $X$ definable over $A$ is $p$-internal if for every tuple $a$
of $X$, $tp(a/A)$ is $p$-internal.
If we replace $dcl$ by $acl$ above we say that $q$ (or $X$) is almost $p$-internal.
\end{definition}

\begin{definition}
Let $p$ be a (possibly partial) type over $A$ and $q=tp(a/B)$ a type. 
The $p$-weight of $q$, denoted by $w_p(q)$, is 
the largest integer $n$ such that there are $C\supset A\cup B$, a tuple $a_1,\ldots,a_n$ of realizations of $p$ 
which are independent over $C$, and a realization $b$ of $q$ such that $(a_1,\ldots,a_n)\dnfo_A C$, $b\dnfo_B C$ 
and $a_i \dfo_C b$ for every $i=1,\ldots,n$. 
If $p$ is the partial type $x=x$ we say weight instead of $p$-weight and it is denoted by $w(q)$.
\end{definition}

\begin{definition}\label{prt11}
Let $A$, $B$ and $C$ be sets. We say that $A$ dominates $B$ over $C$ if for every set $D$, $D \dnfo_C A$ implies
$D \dnfo_C B$.
Let $p,q$ be two types. We say that $p$ dominates $q$
if there is a set $C$ containing the domains of $p$ and $q$ and realizations $a$ and $b$ of non-forking extensions of
$p$ and $q$ to $C$ respectively, such that $a$ dominates $b$ over $C$.
We say that $p$ and $q$ are equidominant if $p$ dominates $q$ and $q$ dominates $p$.
\end{definition}

A type is said to be regular if it is orthogonal to all its forking extensions.
For regular types, equidominance is an equivalence relation and non-orthogonality is transitive (\cite{wag}, section 5.2).

Let $p$ and $q$ be two complete types. We say that $q$ is hereditarily orthogonal to $p$ if every extension of $q$
is orthogonal to $p$.

\begin{definition}\label{prt8}
\begin{enumerate}
\item Let $p$ be a type and $A$ a set. The $p$-closure of $A$, $cl_p(A)$ is the set of all $a$ such that $tp(a/A)$ is hereditarily orthogonal to $p$.
\item A type $p$ is called locally modular if for any $A$ containing the domain of $p$,
 and any tuples $a$ and $b$ of realizations of $p$,
 we have $a \dnfo_C b$ where
$C=cl_p(Aa) \cap cl_p(Ab)$
\end{enumerate}
\end{definition}

For the dichotomy we talk about 1-basedness, which is related with local modularity.

\begin{definition} 
\begin{enumerate}
\item Let $A \subset M$ and let $S$ be an $(\infty)$-definable set over $A$. We say that $S$ is 1-based if
for every $m,n \in \na $, and $a\in S^m,b \in S^n $, $a$ and $b$ are independent over
$acl(Aa) \cap acl(Ab)$.
\item A type is 1-based if the set of its realizations is 1-based.
\end{enumerate}
\end{definition}

For an important characterization of locally modular types we need to introduce the following definitions.

\begin{definition}\label{prt5}
\begin{enumerate}
\item Let $p \in S(A)$ be regular and let $q$ be a type over a set $B \supset A$. 
We say that $q$ is $p$-simple if there is $C \supset B$ and 
a set $X$ of realizations of $p$ and a realization $a$ of $q$ with $a \dnfo_B C$ such that $tp(a/CX)$ is 
hereditarily orthogonal to $p$. 
\item Let $p$ be a regular type over $A$ and let $q$ be a $p$-simple type. 
We say that $q$ is $p$-semi-regular if it is domination equivalent to a non-zero
power of $p$
\end{enumerate}
\end{definition}

\begin{rem}\label{canbase}
If $a, b$ are tuples of $M$, the canonical base of $a$ over $b$ denoted by $Cb(a/b)$ the smallest algebraically
closed subset of $M$ over which $tp(a/b)$ does not fork.
As $T$ is simple, canonical bases exist (see \cite{simple2} for the details). 
3.3 of \cite{wag} implies that for $B \subset M,$ a tuple $a$ of $M$, and  a sequence $a_i$ of realizations of $tp(a/B)$ independent over $B$,
there is  some $m$ for which $Cb(a/B) \sse dcl(a_1 \cdots a_m)$ .

\end{rem}
\begin{prop}
\emph{(\cite{wag},3.5.17)}
\label{prt10}
A type $p$ is locally modular if and only if for any two models $M$ and $N$
with $N \prec M$, and any tuple of realizations $a$ of $p$ over $M$ such that $tp(a/N)$
is $p$-semi-regular, $Cb(a/M) \subset cl_p(Na)$.
\end{prop}

\begin{lem}\label{prt011}
Let $A$ be a subset of $M$, $a$ a tuple and $b\in acl(Aa)$.
Let $q=tp(a,b/A)$ be a regular type and let  $p=tp(b/A)$. Then $p$ is locally modular if and only if $q$ is locally modular.
\end{lem}
{\it Proof:}\\

As being hereditarily orthogonal to $p$ is the same
as being hereditarily orthogonal to $q$, by definition, for any set $B$,
$cl_p(B)=cl_q(B)$.

Let $(a_1,b_1)$ and $(a_2,b_2)$ be tuples (of tuples) of realizations of $q$. 

{\bf Claim}: $a_i\in cl_p(b_i)$ for $i=1,2$:

{\bf Proof of the claim}:
Suppose that $a_1 \not\in cl_p(b_1)$. Then there exists a non-forking extension
$r$ of $tp(a_1)$ over $B$ such that  $r$ is non-orthogonal to $tp(b_1)$.
Let $c$ and $d$ be realisations of the partial type $r \cup q$. Then
$tp(c,d/B)$ is a forking extension of $q$. The element $c$ can be chosen
dependent over $B$ from a realisation $(c', d')$ of a forking extension of $q$ 
over $B$, and this contradicts the regularity of $q$. The same argument applies
to $a_2$ and $cl_p(b_2)$  and the claim is proved.

The following equation holds:
$$cl_q(a_1,b_1)\cap cl_q(a_2,b_2)=cl_p(b_1)\cap cl_p(b_2)=_{\rm def}C.$$ 

It follows immediately that the local modularity of $q$ implies the local modularity of $p$. 
Conversely, assume that $p$ is locally modular. Then 
$b_1\dnfo_Cb_2$. Let $D=Cb(a_1b_1/acl(Cb_2))$. By \ref{canbase}
$D$ is contained in the algebraic closure of independent realizations of
$tp(a_1b_1/acl(Cb_2))$  this implies that $tp(D/C)$ is almost-internal to the set of conjugates of $tp(a_1/Cb_1)$,
 and is therefore hereditarily orthogonal to $p$.
 Hence $D\subset cl_p(C)=C$, and $a_1b_1\dnfo_Cb_2$. 
A similar reasoning gives that $Cb(a_2b_2/acl(Ca_1b_1))\subset C$.\\
$\Box$

\begin{lem}\label{prt012}
Let  $A=acl(A)$ a subset of  $M$ of $T$ and let $a$ be a tuple in $M$.
Assume that $tp(a/A)$ has SU-rank $\beta+\omega^\alpha=\beta\oplus \omega^\alpha$, with $\beta > \omega^{\alpha}$, and has weight 
$1$. Then there is $b\in acl(Aa)$ such that $SU(b/A)=\omega^\alpha$.

\end{lem}

{\it Proof:}\\

By \ref{wgeq}, there is some $C=acl(C)\supset A$ independent from $a$ 
over $A$ and a tuple $c$ such that $\s(c/C)=\omega^\alpha$, and $c$ and 
$a$ are not independent over $C$. Let $B$ be the algebraic closure of 
$Cb(Cc/acl(Aa))$. Then by \ref{canbase} $B$ is contained in the algebraic closure of 
finitely many (independent over $Aa$) realizations of $tp(Cc/acl(Aa))$, 
say $C_1c_1,\ldots,C_nc_n$. Let $D=acl(C_1,\ldots,C_n)$. Then $D$ is 
independent from $a$ over $A$, and each $c_i$ is not independent from 
$a$ over $D$. Since $tp(a/A)$ has weight $1$, so does $tp(a/D)$, and 
therefore for each $1<i\leq n$, $c_1$ and $c_i$ are not 
independent over $D$. Thus $\s(c_i/Dc_1)<\omega^\alpha$, and therefore 
$SU(c_1,\ldots,c_n/D)<\omega^\alpha 2$. As $D$ is independent from $a$ 
over $A$, and $B \subset acl(D,c_1,\cdots,c_n)\cap acl(Aa)$, we get 
$\s(B/A)<\omega^\alpha 2$. 
Since $\s(c/C)=\omega^{\alpha}$ and $\s(c/CB)< \omega^{\alpha}$,
then $\s(B/C)\geq \omega^{\alpha}$, and as $B \dnfo_A C$ we have $\s(B/A)\geq \omega^\alpha$. 
By Lascar's inequalities  we have $\s(a/AB) + \s(B/A) \leq \beta + \omega^{\alpha}$.
As $\s(B/A)\geq \omega^\alpha$ we have that $\s(B/A)=\delta + \omega^{\alpha}$ with $\delta \geq\omega^{\alpha}$ or 
$\delta =0$, 
and $\s(B/A)<\omega^\alpha 2$ implies that  $\delta = 0$.\\
$\Box$

We end this section with two useful results proved in \cite{wag}, section 5.2.

\begin{prop}\label{prt13}
A type in a supersimple theory is equidominant with a finite product of regular types.
\end{prop}

\begin{prop}
\label{orto}
If $p \not\perp q$ and $q \not\perp r$ then there is a conjugate $r'$ of $r$ such that $p \not \perp r'$.
\end{prop}

\section{Difference-differential fields}\label{sec:dcfa}
A difference-differential field is a differential field $(K,D)$ with an automorphism $\si$ of $K$ which commutes with $D$. 

The theory of difference-differential fields of characteristic zero has a model-companion which we denote {\it DCFA} (\cite{rbdcfa1}). 
Before we give an axiomatization of this theory we need to introduce some definitions regarding varieties defined in differential fields.

\begin{definition}\label{I220}
Let $ (K,D) $ be a differential field, and let $V \subset \mathbb{A}^n$
be a variety, let  $F(X)$ be a finite tuple of polynomials over $K$
generating $I(V)$ where $X=(X_1, \cdots , X_n)  $.

\begin{enumerate}

\item  We define the first prolongation of $V$, $\tau_1(V)  $ by the equations:
$$F(X)=0,  J_F(X)Y_1^t +F^D(X)=0  $$
where $Y_1 $ is an $n$-tuple, $F^D$ denotes the tuple of polynomials obtained by applying $D$ 
to the coefficients of each polynomial of $F$,
 and $J_F(X)  $ is the Jacobian matrix of $F$ (i.e. if $F =(F_1, \cdots , F_k) $ then 
$ J_F(X)= (\partial F_i / \partial X_j)_{1 \leq i \leq k,1 \leq j \leq n} $).

\item For $m>1 $, we define the $m$-th prolongation of $V$ by induction on $m$:

Assume that $\tau_{m-1}(V) $ is defined by $F(X)=0$, 
$J_F(X)Y_1^t+F^D(X)=0, \cdots, $\linebreak
$ J_F(X)Y_{m-1}^t+f_{m-1}(X,Y_1,\cdots,Y_{m-2})=0. $
Then $\tau_m(V)  $ is defined by:\\
$$(X,Y_1, \cdots,Y_{m-1}) \in \tau_{m-1}(V)$$ and  
$$ J_F(X)Y_m^t+J_F^D(X)Y_{m-1}^t+J_{f_{m-1}}(X,Y_1,\cdots,Y_{m-2})(Y_1,\cdots,Y_{m-1})^t$$ 
$$+ f_{m-1}^D(X,Y_1,\cdots, Y_{m-2})=0 .$$

\item Let $W \subset \tau_m(V)  $ be a variety. We say that $W$ is in normal form if,
 for every $i \in \{0, \cdots,m-1   \} $,
 whenever $G(X,Y_1, \cdots,Y_i) \in I(W) \cap K [X,Y_1,\cdots,Y_i ]  $ then 

$J_G(X,Y_1,\cdots,Y_i)(Y_1, \cdots, Y_{i+1})^t +G^D(X,Y_1,\cdots ,Y_i) \in I(W). $
\item Let $W \subset \tau_m(V)  $ be a variety in normal form. \\
A point $a$ (in some extension of $K$) is an
$(m,D) $-generic of $W$ over $K$ if\linebreak
$(a,Da, \cdots ,D^ma) $ is a generic of
 $W$ over $K$ and for $i>m  $, we have
 
$\td(D^ia/K(a,\cdots,D^{i-1}a))=\td(D^ma/K(a, \cdots,D^{m-1}a) )  $.
\end{enumerate}
\end{definition}

Now we can give an axiomatization for {\it DCFA}.

\begin{fact}
\emph{\cite{rbdcfa1}.}
\label{I49}
$(K,D,\sigma)$ is a model of  {\it DCFA} if
\begin{enumerate}
\item $(K,D)$ is a differentially closed field.
\item $\sigma$  is an automorphism of $(K,D)$.
\item If $U,V,W$ are varieties such that:
\begin{enumerate}
\item $U \subset V \times V^{\sigma}  $ projects generically onto $V$ and
 $ V^{\sigma}  $.
\item $W \subset \tau_1 (U) $ projects generically onto $U$.
\item $\pi_1(W)^{\sigma}=\pi_2(W)$ \textup{(}we identify $\tau_1(V \times V^{\sigma})$ 
with $\tau_1(V)\times \tau_1(V)^{\si}$\textup{)} and let $\pi_1:\tau_1(V \times V^{\si}) \to \tau_1(V)$
and $\pi_2:\tau_1(V \times V^{\si}) \to \tau_1(V)^{\si}$ be the natural projections\textup{)}. 
\item A $(1,D) $-generic point  of $W$ projects  
onto a $(1,D)   $-generic point of $\pi_1(W)$ and onto a (1,D)-generic point of $\pi_2(W) $. 
\end{enumerate}
Then there is a tuple $a \in V(K) $, such that $(a, \sigma(a)) \in U $ and \\
$(a,Da,\si(a),\si(Da)) \in W$.
\end{enumerate}
\end{fact}

This theory is not complete, but its completions are
easily described. It's properties in fact are in general very similar to the ones of {\it ACFA}: independence is
defined by linear disjointness, the models of {\it DCFA} eliminates imaginaries (moreover, they satisfy the Independence
Theorem over algebraically closed sets), the completions of {\it DCFA} are supersimple and thus types are ranked by
the $\s$-rank; forking is determined by quantifier-free formulas, this implies that {\it DCFA} is
 quantifier-free $\omega$-stable; in a model of {\it DCFA} the difference-differential Zariski topology 
(defined in analogy with Zariski topology
in algebraically closed fields) is Noetherian. All these properties are proved in \cite{rbdcfa1}.

Since {\it DCFA} is quantifier-free $\omega$-stable, we can define canonical bases for quantifier-free
types as in stable theories as in \cite{zopi}. We denote the canonical basis of the quantifier-free type
of $a$ over $K$ as  $Cb(qftp(a/K))$. It does not  coincide with the canonical base of
$tp(a)$ as defined for simple theories. However, as {\it DCFA}
satisfies  the  independence theorem over algebraically closed sets,  $Cb(p)$ will be contained in
$acl(Cb(qftp(a/K))$.

\section{Arc spaces in difference differential fields}\label{sec:arcs}
 
First we define algebraic arc spaces, for the formal definition we refer to \cite{arcs}.

Let $K$ be a field, and $K^{(m)}$ the $K$-algebra
$K[\epsilon]/(\epsilon^{m+1})$. Then, identifying $K^{(m)}$ with
$K\cdot 1\oplus K\cdot \epsilon \ldots \oplus K\cdot \epsilon^m$, we see
that the $K$-algebra $K^{(m)}$ is quantifier-free interpretable in $K$,
if we encode elements of $K^{(m)}$ by $(m+1)$-tuples of $K$.

Let $V\subset \mathbb{A}^\ell$ be a variety defined over $K$.
For $m\in \na$, we
consider
the set $V(K^{(m)})$ of $K^{(m)}$-rational points of $V$.

Using the quantifier-free interpretation of $K^{(m)}$ in $K$, we may identify $V(K^{(m)})$ with a subvariety $\mathcal{A}_mV(K)$ of
$\mathbb{A}^{(m+1)\ell}(K)$.  The variety $\mathcal{A}_mV$ is called the $m$-th arc
 bundle of
$V$. More precisely, if $f_1,\ldots,f_k\in
K[X_1,\ldots,X_\ell]$ generate the ideal $I(V)$, then the ideal of
$\mathcal{A}_mV$ is generated by the polynomials $f_{j,t}\in K[X_{i,t}: 1\leq
i\leq\ell, 0\leq t\leq m]$, $1\leq j\leq k, 0\leq
t\leq m$, which are defined by the identity
$$f_j((\sum_{t=0}^m x_{i,t}\epsilon^t)_{1\leq i\leq l})=\sum_{t=0}^{m}
f_{j,t}((x_{i,t})_{1\leq i\leq l})\epsilon^t$$
modulo $(\epsilon^{m+1})$

If $m=1$ the polynomials defining $\mathcal{A}_1V$ are
determined by $f_j((x_{i,0} +x_{i,1}\epsilon)_{1\leq i\leq l})=
f_{j,t}((x_{i,0})_{1\leq i\leq l})+
f_{j,t}((x_{i,1})_{1\leq i\leq l})\epsilon$ modulo 
$(\epsilon)^2$ and we identify  $\mathcal{A}_1V$ with the tangent bundle $T(V)$.

If $r>m$,  the natural map
$K^{(r)}\to K^{(m)}$  induces a map $V(K^{(r)})\to V(K^{(m)})$, which
in turn induces a morphism $\rho_{r,m}: \mathcal{A}_rV\to \mathcal{A}_mV$. 

In general, given a morphism of varieties $f:U\to V$  defined over $K$, the
natural morphism $U(K^{(m)})\to V(K^{(m)})$ induced by $f$ gives
rise to a morphism $\mathcal{A}_mf:\mathcal{A}_mU\to \mathcal{A}_mV$.

Let us  write $\rho_m$ for $\rho_{m,0}$. For $a\in V(K)$ the $m$-th
arc space of $V$ at $a$, $\mathcal{A}_mV_a$ is the fiber of $\rho_m$ over $a$.

We follow \cite{arcs} to prove our theorem. Our propositions
\ref{arc1}, \ref{arc2}, \ref{arc3}, \ref{arc4}, \ref{arc5}, \ref{arc7}, \ref{arc9}, \ref{arc10}, \ref{arc11}, \ref{arc12}, \ref{arc14}, \ref{arc15} and \ref{arc16} 
are difference-differential versions of 2.9, 2.10, 2.11, 2.12,  3.1, 3.6, 3.9, 3.10,
3.11, 3.12, 3.15 and 3.16 of \cite{arcs}. Some of the  proofs are very similar, we give the 
details where the proofs are different or in case they give some useful insights.

 The following three lemmas are crucial, and allow us to characterise varieties by their arc spaces.

\begin{lem}
\emph{\cite{arcs}.}
\label{arc01}
Let $U,V$ be two algebraic varieties, and let $f:U \to V$ be a morphism, all defined over $K$.
Let $m \in \na$ and $a_m \in \mathcal{A}_mU(K)$ be such that  $\bar{a}=\rho_m(a)$
and $\bar{f(a)}=\rho_m(f(a))$ are non-singular. Let $U'$ be the fiber of
$\rho_{m+1,m}:\mathcal{A}_{m+1}U \to \mathcal{A}_mU$ over $a$ and $V'$ the fiber of
$\rho_{m+1,m}:\mathcal{A}_{m+1}V \to \mathcal{A}_mV$ over $\mathcal{A}_mf(a)$.
Let $\bar{a}=\rho_m(a)$. Then there are biregular maps $\varphi_U:U' \to T(U)_{\bar{a}}$  and
$\varphi_V:V' \to T(V)_{f(\bar{a})}$  such that the following diagram is commutative:

$$\xymatrix{
  & U' \ar[d]^{\varphi_U} \ar[r]^{\mathcal{A}_m(f)}
                 & V' \ar[d]^{\varphi_V}       \\
  & T(U)_{\bar a} \ar[r]^{df_{\bar a}}   & T(V)_{f({\bar a})}                }$$
\end{lem}

\begin{lem}
\emph{\cite{arcs}.}
\label{arc02}
Let $U,V$ be algebraic varieties defined over $K$, and let $f:U \to V$ be a dominant map defined over $K$.
Let $a \in U(K)$ be non-singular such that $f(a)$ is non-singular and the rank of $df_a$  equals $\dim V$. Then for every $m \in\na$
the map $\mathcal{A}_m(f):\mathcal{A}_mU_a(K) \to \mathcal{A}_mV_{f(a)}(K)$ is surjective.
\end{lem}

\begin{lem}
\emph{\cite{arcs}.}
\label{arc03}
Let $U,V,W$ be algebraic varieties defined over $K$ such that $U,V \subset W$. Let $a\in U(K) \cap V(K)$ be non-singular.
 Then
$U=V$ if and only if $\mathcal{A}_mU_a(K)=\mathcal{A}_mV_a(K)$ for all $m \in \na$.
\end{lem}

Let $(\mathcal{U},\si,D)$ be a saturated model of {\it DCFA}, let $K$ be a
difference-differential subfield of ${\mathcal U}$. 

Let $V$ be an algebraic variety of the affine space of dimension $k$
over $\mathcal U$.

Let $\nabla_m:V \to \tau_m(V)$ be defined by $x \mapsto (x,Dx,\cdots,D^mx)$ and
let $\pi_{l,m}:\tau_l(V) \to \tau_m(V)$ be the natural projection for $l \geq m$.
$S_m(V)$ will
denote the Zariski closure of $\{(x, \cdots, \si^m(x)): x \in V\}$.
Let $q_m:V \to S_m(V)$ be defined by $x \mapsto (x, \cdots, \si^m(x))$ and let
$p_{l,m}:S_l(V) \to S_m(V)$ be the natural projections for $l \geq m$.\\

We now define a notion of difference-differential prolongation. 
Let $\Phi_m(V)=\tau_m(S_m(V))$, let $\psi_m:V \to \Phi_m(V)$ be such that
$x \mapsto \nabla_m(q_m(x))$ and for $l \geq m$ let
$t_{l,m}:\Phi_l(V) \to \Phi_m(V)$ be defined by $t_{l,m}=\pi_{l,m} \circ p_{l,m}$.
Let us denote $\pi_l=\pi_{l,0}$, $p_l=p_{l,0}$, $t_l=t_{l,0}$,
$\Phi(V)=\Phi^1(V)=\Phi_1(V)$ and $\Phi^{m+1}(V)=\Phi(\Phi^m(V))$.

We define $\psi=\psi^1=\psi_1:V \to \Phi(V)$ and
$\psi^{m+1}(V)=\psi(\psi^m):V \to \Phi^{m+1}(V)$.\\

We extend $\si$ and $D$ to $K^{(m)}$ by defining $\si(\epsilon)=\epsilon$
and $D \epsilon=0$. Since $\si$ and $D$ commute, $\nabla_m q_m(x)$ is a permutation of $q_m \nabla_m x$. For $x=(x_0, \cdots, x_m) \in K^{(m)}$ we identify it with $ \sum_{i=0}^m x_i \epsilon^i$, we can identify $\Phi_m(x)$ with 
$\sum_{i=0}^m \Phi_m(x_i) \epsilon^i$
 Then
we can identify $\mathcal{A}_r(S_m(V))(K)$ with $S_m(\mathcal{A}_r(V))(K)$.
We can thus assume that $\mathcal{A}_r(\Phi_m(V))(K) = \Phi_m(\mathcal{A}_m(V))(K)$.
All this is done with greater generality in \cite{moosa_scanlon_2009} and \cite{moosa_scanlon_2010}
\\

Let $X$ be a $(\si,D)$-variety given as a $(\si,D)$-closed subset of an
algebraic variety $\bar{X}$. We define $\Phi_m(X)$ as the Zariski closure
of $\psi_m(X)$ in $\Phi_m(\bar{X})$. Thus $X$ is determined by the prolongation sequence
$\{ t_{l,m}:\Phi_l(X) \to \Phi_m(X)):l \geq m \}$, since
$X(\mathcal{U})=\{a \in {\bar X}({\mathcal U}): \psi_l(a) \in \Phi_l(X) \forall l \}$.
We call this sequence the prolongation sequence of $X$. \\

\begin{prop}\label{arc04}
Let $\{ V_l \subset \Phi_l(\bar{V}):l \geq 0\}$ be a sequence of algebraic varieties
and $\{t_{m,l}:V_m \to V_l, m\geq l\}$ a sequence of morphisms such that:
\begin{enumerate}
\item $t_{l+1,l} \upharpoonright V_{l+1} \to V_l$ is dominant.
\item After embedding $\Phi_l(\bar{V})$ in $\Phi^l(\bar{V})$ and
$\Phi_{l+1}(\bar{V})$ in $\Phi^{l+1}(\bar{V})$,
\begin{enumerate}
\item $V_{l+1}$ is a
subvariety of $\Phi(V_l)$.
\item Let $\pi'_1:\Phi(V_l)\to \tau(V_l)$ and $\pi'_2:\Phi(V_l)\to
\tau(V_l^\sigma)$ be the projections induced by $\Phi(V_l)\subset
\tau(V_l)\times \tau(V_l^\sigma)$; then $\pi'_1(V_{l+1})^\sigma$ and
$\pi'_2(V_{l+1})$ have the same Zariski closure.
\end{enumerate}
\end{enumerate}
Then there is a (unique) $(\si,D)$-variety $V$ with prolongation sequence $\{t_{m,l}:V_m \to V_l, m\geq l\}$ .
\end{prop}
{\it Proof:}\\

For each $l$, as the maps $\pi_{m,j}$ are dominant, the system $\{p_{m,l}(V_m), \pi_{m,j}:m >j \geq l\}$
defines a differential subvariety  $W_l$ of $\bar{V} \times \cdots \times \bar{V}^{\si^l}$.

Condition (1). implies that for $m$ sufficiently large, an $(m,D)$-generic
of $p_{m,l+1}(V)$ is sent by $p_{l+1,l}$ to an $(m,D)$-generic of 
$p_{m,l}(V)$. Hence,
a $D$-generic of $W_{l+1}$ is sent by $p_{l+1,l}$ to a $D$-generic of 
$W_l$ (that is, a generic in the sense of {\it DCF}).

By conditions (2) (b) and (1), the  map $t'_{l+1,l}:V_{l+1}\to V_l^\sigma$ 
induced
by $\Phi(V_l)\to V_l^\sigma$ is dominant. Hence, considering the 
natural
projection
 $p'_{l+1,l}:S_{l+1}(\bar
V)\to S_l(\bar V)^\sigma$, and reasoning as above, we obtain that
$p'_{l+1,l}$ sends a $D$-generic of $W_{l+1}$ to a $D$-generic of
$W_l^\sigma$. 
Hence by \ref{I49}, for every $l$ there is $a$ such that $\psi_l(a)$ is a generic of $V_l$ over $K$. 
By saturation, there is $a$ such that for all $l$, $\psi_l(a)$ is a generic of $V_l$. Then
$\{t_{m,l}:V_m \to V_l, m\geq l\}$  is the prolongation sequence of the $(\si,D)$-locus of $a$ over $K$.\\
$\Box$

We define non-singular points in analogy with the  corresponding notion in \cite{arcs}.
\begin{definition}
Let $X$ be a $(\si,D)$-subvariety of the algebraic variety $\bar{X}$.
We say that a point $a \in X$ is non-singular if, for all $l$, $\psi_l(a)$ is
a non-singular point of $\Phi_l(X)$, the maps $dt_{l+1,l}$ and $dt'_{l+1,l}$ at $\psi_{l+1}(a)$ have rank
equal to $\dim X_l$ and the maps $d \pi_1'$ and $d\pi_2'$ (as defined above) at $\psi_{l+1}(a)$ have rank equal to 
the dimension of the Zariski closure of $\pi_1'(\Phi_{l+1}(X))$.
\end{definition}

The next proposition is essentially 4.7 of \cite{moosa_scanlon_2010}.
\begin{prop}\label{arc1}
Let $(K,\si,D)$ be a model of {\it DCFA}. Let $V$ be a $(\si,D)$-variety
given as a closed subvariety of an algebraic variety $\bar{V}$.
Let $m \in \na$ and $a\in V(K)$ a non-singular point.
Then $\{\mathcal{A}_m(t_{r,s}): \mathcal{A}_m\Phi_r(V)_{\psi_r(a)} \to \mathcal{A}_m\Phi_s(V)_{\psi_s(a)}, r \geq s \}$
 form the
 $(\si,D)$-prolongation sequence of
a $(\si,D)$-subvariety of $\mathcal{A}_m\bar{V}_a$.
We define the $m$-th arc space of $V$ at $a$, $\mathcal{A}_mV_a$, to be this subvariety.
We have also that $\Phi_r(\mathcal{A}_mV_a)=\mathcal{A}_m\Phi_r(V)_{\psi_r(a)}$ for all $r$.
\end{prop}
{\it Proof:}\\

As $\mathcal{A}_m\Phi_r(\bar{V}) = \Phi_r(\mathcal{A}_m{\bar V})$, we look at $\mathcal{A}_m\Phi_r(V)_{\psi_r(a)}$ as an
algebraic subvariety of $\Phi_r(\mathcal{A}_m \bar{V})_{\psi_r(a)}$.
We have that $\Phi_{r+1}(V) \subset \Phi(\Phi_r(V))$ for all $r$. Since
$\mathcal{A}$ preserves inclusion we have
$\mathcal{A}_m\Phi_{r+1}(V)_{\psi_{r+1}(a)} \subset \mathcal{A}_m\Phi(\Phi_r(V)_{\psi_{r}(a)})=\Phi(\mathcal
{A}_m\Phi_r(V)_{\psi_{r}(a)})$. This proves conditions 1. and 2.(a) of \ref{arc04}.

Moreover, the maps $t_{r,s}:\Phi_r(V) \to \Phi_s(V)$
are dominant, and as $a$ is non-singular,  by \ref{arc02}, the maps
$\mathcal{A}(t_{r,s}):\mathcal{A}_m\Phi_r(V)_{\psi_r(a)} \to \mathcal{A}_m\Phi_s(V)_{\psi_s(a)}$,
are dominant.
Applying $\mathcal{A}_m$ to the dominant maps $\pi_1':\Phi_{r+1}(V) \to \tau(\Phi_r(V))$ 
and $\pi_2':\Phi_{r+1}(V) \to \tau(\Phi_r(V))^{\si}$,
using the hypothesis on $a$ and \ref{arc02}, we get
$$\mathcal{A}_m\pi_1'(\mathcal{A}_m(\Phi_{r+1}(V)_{\psi_{r+1}(a)}))=
\mathcal{A}_m(\pi_1'(\Phi_{r+1}(V))_{\pi_1'(\psi_{r+1}(a))})$$
and
$$\mathcal{A}_m\pi_2'(\mathcal{A}_m(\Phi_{r+1}(V)_{\psi_{r+1}(a)}))=
\mathcal{A}_m(\pi_2'(\Phi_{r+1}(V))_{\pi_2'(\psi_{r+1}(a))})$$
and since $\pi_1'(\Phi_{r+1}(V))^{\si}$ and $\pi_2'(\Phi_{r+1}(V))$ have the same Zariski closure,
and $\si(\pi_1'\psi_{r+1}(a))=\pi_2'\psi_{r+1}(a)$ we get condition $2(b)$.
Hence $\{
\mathcal{A}_m(t_{r,s}):A_m\Phi_r(V)_{\psi_r(a)} \to  A_m\Phi_s(V)_{\psi_s(a)}, r \geq s  \}$ is the $(\si,D)$-prolongation
sequence of a $(\si,D)$-subvariety $W$ of $\mathcal{A}_m\bar{V}_a$, where
$W(K)=\{x \in \mathcal{A}_m\bar{V}_a(K): \psi_r(x) \in \mathcal{A}_m\Phi_r(V)_{\psi_r(a)}(K), 
r \geq 0\}$ and $\mathcal{A}_m\Phi_r(V)_{\psi(a)}=\Phi_r(W)$ for all $r$. We define then
$\mathcal{A}_mV_a=W$.\\
$\Box$

\begin{lem}\label{arc2}
Let $U,V$ be two $(\si,D)$-subvarieties of an algebraic variety $\bar{V}$.
Let $a \in U(K) \cap V(K)$ be a non-singular point of $U$. Then
$U=V$ if and only if $\mathcal{A}_m\Phi_l(U)_{\psi_l(a)}=\mathcal{A}_m\Phi_l(V)_{\psi_l(a)} $ for all $m,l$.
\end{lem}
{\it Proof:}\\

If $\mathcal{A}_mU_a(K)=\mathcal{A}_mV_a(K)$ for all $m$, then
$\Phi_r(\mathcal{A}_mU_a)(K)=\Phi_r(\mathcal{A}_mV_a)(K)$. Thus, by \ref{arc1},
$\mathcal{A}_m\Phi_r(U)_{\psi_r(a)}(K)=\mathcal{A}_m\Phi_r(V)_{\psi_r(a)}(K)$.
Hence, for all $r$ and $m$, we have
$\mathcal{A}_m\Phi_r(U)_{\psi_r(a)}=\mathcal{A}_m\Phi_r(V)_{\psi_r(a)}$.
Lemma \ref{arc03} implies that $U$ and $V$ have the same $(\si,D)$-prolongation sequence. Hence $U=V$.\\
$\Box$

\begin{definition}
Let $V$ be a $(\si, D)$-variety and $a$ a non-singular point of $V$. We define the $(\si,D)$-tangent space
 $T_{\si,D}(V)_a$ of $V$ at $a$
as follows:

Let $P_r$ be a finite tuple of polynomials generating $I(\Phi_r(V)_{\psi_r(a)})$. Then
$T_{\si,D}(V)_a$ is defined by the equations $J_{P_r}(\psi_r(a)) \cdot (\psi_r(Y))= 0$.
In other words, the prolongation sequence of $T_{\si,D}(V)_a$ is 
$dt_{l,r}: T(\Phi_l(V))_{\psi_l(a)} \to T(\Phi_r(V))_{\psi_r(a)},l \geq r\}$,
where $T$ denotes the usual tangent bundle and $t_{l,r}$ the natural projection 
$\Phi_l(V)_{\psi_l(a)} \to \Phi_r(V)_{\psi_r(a)}$.
\end{definition}

\begin{rem}
Let $a$ be a non-singular point of the $(\si,D)$-variety $V$. Then $T_{\si,D}(V)_a$ is a subgroup
of $\gr_a^n(K)$, and by the same arguments as above, its prolongation sequence is
$(d(t_{l,r})_{\psi_l(a)}:T(\Phi_l(V))_{\psi_l(a)} \to T(\Phi_r(V))_{\psi_r(a)})_{l \geq r}$.
\end{rem}
\begin{lem}\label{arc3}
Let $V$ be a $(\si,D)$-variety in $\mathbb{A}^l$ and $a$ a non-singular point of $V$. Then $\mathcal{A}_1V_a$
is isomorphic to $T(V)_a$. Let  $\bar{V}$ be the Zariski closure of 
$V(\mathcal{U})$ in $\mathbb{A}^l$ and $m \in \na$; 
then the map given by lemma \ref{arc01} which identifies the fibers of
$\mathcal{A}_{m+1}\bar{V}_a \to \mathcal{A}_m\bar{V}_a$ with $T(\bar{V})_a$ restricts
to an isomorphism of the fibers of
$\mathcal{A}_{m+1}V_a \to \mathcal{A}_mV_a$ with $T(V)_a$.
\end{lem}

{\it Proof:}\\

As remarked before, we can assume that $\mathcal{A}_1\bar{V} =T(\bar{V})$.
Let $b \in T(\bar{V})_a(K)$. By definition
$(a,b) \in \mathcal{A}_1V(\mathcal{U})$ if and only if $\psi_r(a,b) \in T(\Phi_r(V))(K)$
 for all $r$. As  $T(\Phi_r(\bar{V}))= \Phi_r(T(\bar{V}))$, $T(\Phi_r(V))_{\psi_r(a)}$ is an algebraic subvariety of
 $\Phi_r(T(\bar{V}))$ and $\psi_r(a,b)=(\psi_r(a),\psi_r(b))$. Hence $b\in \mathcal{A}_1V_a(K)$
 if and only if $b \in T(V)_b$ and the first part of the theorem is proved.

 Now we look at the map given in \ref{arc01}. In particular, if
 $c \in \mathcal{A}_mV_a(K)$ and $r \geq 0$, by \ref{arc1},
 $\psi_r(c) \in \mathcal{A}_m\Phi_r(V)_{\psi_r(a)}$ and the following diagram commutes

$$\xymatrix{
  & (\mathcal{A}_{m+1}{\bar V}_a)_c \ar[d] \ar[r]
                 & (\mathcal{A}_{m+1}\Phi_r({\bar V})_{\psi_r(a)})_{\psi_r(c)} \ar[d]       \\
  & T({\bar V})_{a} \ar[r]   & T(\Phi_r({\bar V}))_{\psi_r(a)}                }$$
where the horizontal arrows are $\psi_r$ and the vertical arrows are the maps given by \ref{arc01} applied to $\bar{V}$
 and $\Phi_r(\bar{V})$. So $(\mathcal{A}_{m+1}V_a)_c$ is
 identified with $T_{\si,D}(V)_a$.\\
$\Box$

\begin{notadefi}

In  \cite{Levin}  the author proved that given a difference-differential
 subfield $F$ of $K$ and $a \in K$ there is a numerical polynomial 
 $P_{a/F}(X) \in {\mathbb Q}[X]$ of degree at most 2, such that for sufficiently large $r \in \na$,
 $P_{a/F}(r)=\td(\psi_r(a)/F)$. We call the degree of $P_{a/F}$ the
 $(\si,D)$-type of $a$ over $F$, and the leading coefficient of $P_{a/F}$
 the typical dimension of $a$ over $F$, it is denoted
 $\dim_{\si,D}(a/F)$. For a $(\si,D)$-variety $V$ defined over $F$ we define
 $P_V=P_{a/F}$ where $a$ is a $(\si,D)$-generic of $V$ over $F$.
 We have that the $(\si,D)$-type of $a$ over $F$ is 2 if and only if $a$ contains an element which
is $(\si,D)$-transcendental over $F$.

Let $F=acl(F) \subset \mathcal{U}$ and
let $a \in \mathcal{U}$ and let $p=tp(a/F)$. We denote by $m(p)$ (or by $m(a/F)$) the $(\si,D)$-type of $a$ over $F$
and we write $\dim_{\si,D}(p)$ for $\dim_{\si,D}(a/F)$.
If $p'$ is a non-forking extension of $p$ then $m(p)=m(p')$ and $\dim_{\si,D}(p) =\dim_{\si,D}(p')$.
If $A$ is an arbitrary subset of $\mathcal U$ we write $m(a/A)$ instead of $m(a/acl(A))$.

If $V$ is a $(\si,D)$-variety over $K$, $m(V)$ denotes the $(\si,D)$-type of $V$. Then, if
$a$ is a $(\si,D)$-generic of $V$, $m(V)=m(qftp(a/F))$.

\end{notadefi}

 \begin{cor}\label{arc4}. 
Let $V$ be a $(\si,D)$-variety in $\mathbb{A}^l$, and $m \in \na$. Then
 for $a \in V(K)$ non-singular, the
$(\si,D)$-type of $V$ and $\mathcal{A}_mV_a$ are equal.
\end{cor}
{\it Proof:}\\

By \ref{arc1} $\Phi_r(\mathcal{A}_mV_a)=\mathcal{A}_m\Phi_r(V)_{\psi_r(a)}$. But if $b$ is a non-singular
point of a variety $U$, then by \ref{arc3}, we have that $\dim(\mathcal{A}_mU_b)=m \dim(U)$.\\
$\Box$

\begin{rem}\label{arc41}
By \ref{arc3}, if $m=1$ and for $a \in V(K)$ non-singular, we have $P_V=P_{T(V)_a}$.
\end{rem}

\begin{lem}\label{arc5} Let $F=acl(F)$. Then
\begin{enumerate}
\item $m(a,b/F)=max\{m(a/F),m(b/F) \}$.
\item If $m(a/F)=m(b/F)$ then $\dim_{\si,D}(a,b/F)=\dim_{\si,D}(a/F)+\dim_{\si,D}(b/Fa)$.
\item If $m(a/F)>m(b/F)$ then $\dim_{\si,D}(a,b/F)=\dim_{\si,D}(a/F)$.
\end{enumerate}
\end{lem}
We borrow the next definitions from \cite{arcs}
\begin{definition}\label{arc8}
\begin{enumerate}
\item Let $p$ be a regular type. We say that $p$ is
$(\si,D)$-type minimal if for any type $q$, $p \not\perp q$ implies
$m(q) \geq m(p)$.
\item A $(\si,D)$-variety $V$ is $(\si,D)$-type minimal if for every
proper $(\si,D)$-subvariety $U$, $m(V) < m(U)$.
\end{enumerate}
\end{definition}

\begin{lem}\label{arc6}
Let $p$ be a  type and let $V$ be the $(\si,D)$-locus of $p$ over $K$
If $V$ is $(\si,D)$-type minimal then $p$ is regular and $(\si,D)$-type minimal.
\end{lem}
{\it Proof:}\\

Let $a$ be a realization of a forking extension of $p$ to some
$L=acl(L) \supset K$. Let $b$ realize a non-forking extension of $p$
to $L$. Let $U$ be the $(\si,D)$-locus of $(a,b)$ over $L$. Then the projection
on the second coordinate: $U \to V$ is dominant, thus $m(a,b/L) \geq m(V)$.
Now if $a \dfo_L b$, then the $(\si,D)$-locus of $b$ over
$acl(La)$  is a proper subvariety of $V$ and therefore, by \ref{arc5},  $m(b/La)<m(V)$; from $m(a/L)<m(V)$, we deduce
$m(a,b/L)<m(V)$ which is impossible.\\
$\Box$

\begin{lem}\label{arc7}
If $p$ is a type over $K$, there is a finite sequence of regular types $p_1, \cdots, p_k$ such that $m(p) \geq m(p_i)$ for all $i$ and $p$
is domination-equivalent to $p_1 \times \cdots \times p_k$.
\end{lem}
{\it Proof:}\\

By \ref{prt13} it suffices to show that given a regular type $q$, such that $p \not\perp q$, there
is a regular type $r$ such that $q \not\perp r$ and $m(r) \leq m(p)$.
Let $a$ be a realization of a non-forking extension of $p$ to some $L$ and let
$b$ be a realization of a non-forking extension of $q$ to $L$ such that
$a \dfo_L b$. Let $c=Cb(tp(a/L,b))$. Thus $c \not\in acl(L)$ and
$c \in acl(Lb)$. So $r=tp(c/L)$ is  non-orthogonal to $q$ and regular (because $c \in acl(Lb)$ and regularity is preserved by algebraicity) .
On the other hand, by \ref{canbase}, there are $a_1, \cdots,a_l$ realizations of $p$ such that
$c \in dcl(La_1 \cdots a_l)$. Then, by \ref{arc5}, $m(r) \leq m(q)$.\\
$\Box$

\begin{rem}\label{remrank}
By \ref{prt012}, in the proof above we can suppose that $\s(c/L)= \omega^i$ for $i \in \{0,1,2 \}$. Thus, given a type $p$  over $K$. there is a type $q$ such that
it is $(\si, D)$-type minimal and has $\s$-rank $\omega^i$.
\end{rem}

\begin{lem}\label{arc9}
Let $G$ be a $(\si,D)$-vector group (that is, a $(\si,D)$-variety which is a subgroup of
$\gr_a^k$ for some $k$). Then $T_{\si,D}(G)_0$ is definably isomorphic to $G$.
Moreover, if $H$ is a $(\si,D)$-subgroup of $G$, then the
restriction of this isomorphism to $H$ is an isomorphism between $H$ and $T_{\si,D}(H)_0$.
\end{lem}
{\it Proof:}\\

Suppose that $G$ is a $(\si,D)$-subgroup of $\gr_a^k$. For each $r \in \na$, $\Phi_r(G)$ is a subgroup
of $\Phi_r(\gr_a^k)=\gr_a^{k(r+1)^2}$. 
Since $\Phi_r(G)$ is an algebraic subgroup of $\gr_a^{k(r+1)^2}$ , its defining ideal is generated by
homogeneous linear  polynomials, and thus its tangent space at 0 is defined by the same polynomials.
Then the map $z_r:\Phi_r(\gr_a^k) \to T(\Phi_r(\gr_a^k))$ defined by
$x \mapsto(0,x)$ identifies $\Phi_r(\gr_a^k)$ and $T(\Phi_r(\gr_a^k))_0$ and
 it restricts to an isomorphism $\Phi_r(G) \to T(\Phi_r(G))_0$. Hence $(z_r:r \geq 0)$
identifies the prolongation sequence of $G$ and the prolongation sequence of $T_{\si,D}(G)_0$.

For the moreover part, it suffices to note that, by our construction above, the restriction
of $z_r$ to $\Phi_r(H)$ is an isomorphism between  $\Phi_r(H)$ and $T(\Phi_r(H))_0$. \\
$\Box$

\begin{cor}\label{arc10}
Let $G$ be a $(\si,D)$-subgroup of $\gr_a^k$. Suppose that for every proper definable subgroup $H$ of $G$,
$m(H)<m(G)$. Then $m(V)<m(G)$ for any proper $(\si,D)$-subvariety of $G$. In particular the generic type
of $G$ is regular.
\end{cor}
{\it Proof:}\\

Let $V$ be a  $(\si,D)$-type minimal  $(\si,D)$-subvariety of $G$ such that $m(V)=m(G)$. After possibly replacing
$V$ by a translate we may assume that $0 \in V$ and $0$ is non-singular. By Remark \ref{arc41}, $m(T(V)_0)=m(V)=m(G)$.
Since  $T(V)_0$ is a subgroup of $T(G)_0 \simeq G$,we obtain $T(V)_0=T(G)_0$.
By \ref{arc4}, $P_V=P_{T(V)_0}=P_{T(G)_0}=P_G$. Hence $V=G$. By \ref{arc6}, the generic type of $G$
is regular.\\
$\Box$

\begin{lem}\label{arc11}
Let $a,c$ be tuples of $\mathcal{U}$. Let $V$ be the $(\si,D)$-locus of $a$ over $K$. Assume that
$c=Cb(qftp(a/acl(Kc)))$. Then there is $m \in \na$ and a tuple $d$ in $\mathcal{A}_mV_a$ such that
$c \in K(a,d)_{\si,D}$.
\end{lem}
{\it Proof:}\\

Let $U$ be the $(\si,D)$-locus of $a$ over $acl(Kc)$. Then $\mathcal{A}_mU_a \subset \mathcal{A}_mV_a$.
 As {\it DCFA} eliminates imaginaries every definable set has a canonical parameter. 
Then $c$ is interdefinable with the canonical parameter of $U$ which,
by \ref{arc2}, is interdefinable over $K(a)_{\si,D}$ with the sequence of the canonical parameters of $\mathcal{A}_mU_a$ over
$K(a)_{\si,D}$. By quantifier-free stability $\mathcal{A}_mU_a$ is defined with parameters from 
$\mathcal{A}_mU_a$.\\
$\Box$

\begin{lem}\label{arc12}
Let $(K,\si,D)$ be a submodel of $(\mathcal{U},\si,D)$. Let $V$ be a $(\si,D)$-variety defined over $K$ and
let $a \in V(\mathcal{U})$ be a non-singular point. Let $b \in \mathcal{A}_mV_a$. Then there are $b_1,b_2, \cdots,b_m=b$,
such that
$b_i \in acl(Ka,b)$ and each $b_i$ is in some $K\cup \{a,b_{i-1}\}$-definable principal homogeneous space for
$T(V)_a$.
\end{lem}
{\it Proof:}\\

By \ref{arc01} and \ref{arc3} each fiber $\rho_{i+1,i}:\mathcal{A}_{i+1}V_a \to \mathcal{A}_iV_a$ is
a principal homogeneous space for $T(V)_a$. Then set $b_i=\rho_{m,i}(b)$.\\
$\Box$

The following lemmas will state the connections between regular non locally modular types and vector groups.

\begin{lem}\label{arc13}
Let $(K,\si,D)$ be a submodel of $(\mathcal{U},\si,D)$. Let $p$ be a $(\si,D)$-type minimal
regular type over $K$ such that $m(p)=d$. If $p$ is not locally modular, 
then there are a vector group $G$ and a type
$q$  such that:
\begin{enumerate}
\item $m(q)=m(G)=d$.
\item $(x \in G) \in q$.
\item $p \not\perp q$.
\end{enumerate}
\end{lem}
{\it Proof:}\\

By \ref{remrank}  we may assume that $\s(p)=\omega^j$ where $j\in\{0,1,2\}$.
By \ref{prt10}, enlarging $K$ if necessary, there are tuples 
$a$ and $c$, with $a$ a tuple of realizations of $p$, $tp(c/K)$ $p$-internal, 
$c=Cb(a/acl(Kc))$, $tp(a/Kc)$ $p$-semi-regular and $c\notin cl_p(Ka)$.
Let $V$ be the locus of $a$ over $K$.

By \ref{arc11}, there is a $k$-tuple $d$ in $\mathcal{A}_mV_a(\mathcal{U})$ such that $c\in
acl(K,a,d)$. For $i=1,\ldots,m$ let $d_i=\rho_{m,i}(d)$. Then for each
$i$, $d_i$ is in some $K(ad_{i-1})$-definable $T(V)_a^k$-principal
homogeneous space. 

Let $m=w_p(c/Ka)$. This means that for any $L=acl(L) \subset \mathcal{U}$
such that $L \dnfo_K c$, given a tuple $(g_1, \cdots, g_m)$
realizing $p^{(m)}$ we have that $g_i \dfo_L c$ for all $i$ if
and only if $g \subset cl_p(Lc)$.

As $c\in acl(K,a,d)$, $c \not\in cl_p(K,a)$, $tp(c/K)$ is $p$-internal;
there is $j \in \{1, \cdots, k\}$ such that $w_p(c/Kad_{j-1})=m$ and
$w_p(c/Kad_j) \leq m-1$.
 Let $L=acl(L) \subset \mathcal{U}$ contain $Kad_{j-1}$,
such that $L \dnfo_K C$, and $(g_1, \cdots, g_m)$ 
realizing  $p^{(m)}$ such that $g_i \dfo_L c$ for all $i$. 
Since  $w_p(c/Kad_{j-1})>w_p(c/Kad_j)$, either there  is $g_k$ such that
$g_k \dnfo_{Ld_j} c$, or $tp(g_k/Ld_j)$ forks over $L$.
In both cases,  $d_j$ and $g$ are dependent
over $L$. Hence
$tp(d_j/Kad_{j-1})\not\perp p$. 

Let $q=tp(d_j/Kad_{j-1})$.
Then we have $m(q)=m(H) \leq m(T(V)_a)=m(p)$, hence $m(p)=m(q)$.\\
$\Box$

\begin{lem}\label{arc14}
Let $p$ be a  regular $(\si,D)$-type minimal type. If there are a $(\si,D)$-vector group $G$
and a type $q$ that satisfy the conclusions of  \ref{arc13}, then there exists a $(\si,D)$-vector group
whose generic type is regular, $(\si,D)$-type minimal and non-orthogonal to $p$.
\end{lem}
{\it Proof:}\\

We order the triplets $ord(G)=\{m(G),\dim_{\si,D}(G),\s(G)\}$ with the lexicographical
order. We proceed by induction on $ord(G)$.

{\bf Claim}:
We may assume that if $H$ is a proper quantifier-free connected,
 quantifier-free
definable subgroup of $G$,  then $m(H)<m(G)$.

{\bf Proof of the claim}:
Suppose that $m(H)=m(G)$. Let $\mu:G \to G/H$ be the quotient map. By \ref{arc5}, $ord(G)>ord(G/H)$.
If we replace $q$ by a non-forking extension of $q$ we may assume that $H$ is defined over the domain $A$ of $q$.
Let $a$ be a realization of $q$ with $tp(a/A) \not\perp p$. As
$q \not\perp p$, we have either $p \not\perp q_0=qftp(\mu(a)/A)$ or
$p \not\perp q'=qftp(a/A\mu(a))$.
If $p \not\perp q_0$  then $m(p) \leq m(q_0)$ by \ref{arc8},  and since $(x \in G/H) \in q_0$,
$m(q_0)\leq m(G/H)\leq m(G)=m(p)$. So $m(q_0)=m(p)$ and we apply the induction hypothesis to $p,q_0$ and $G/H$.
If $p\not\perp q'$, let $b$ be a
realization of $qftp(a/A\mu(a))$ such that $b \dnfo_{A \mu(a)} a$.
Then  $a-b \in H$ and $p \not\perp q''=qftp(a-b/Ab)$ and the same
argument applies.\\

By \ref{arc10} and as $q$ is realized in $G$ and
$m(p)=m(q)=m(G)$, $q$ is a generic of $G$, and is regular and $(\si,D)$-type minimal.
$\Box$

\begin{cor}\label{arc15}
Let $p$ be regular non locally modular  type. Then there is a $(\si,D)$-vector group $G$ whose generic type
is $(\si,D)$-type minimal and non-orthogonal to $p$.
\end{cor}
{\it Proof:}\\

By \ref{arc7} there is a regular type $q$ of minimal $(\si,D)$-type which is non-orthogonal to $p$.
By \ref{arc13}, $q$ satisfies the hypothesis of \ref{arc14}, then there is a $(\si,D)$-vector group $G$ whose generic type
$r$ is nonorthogonal to $q$; again by \ref{orto}, then there is such an $r$ which is non-orthogonal to $p$.\\
$\Box$

\begin{lem}\label{arc16}
Let $G$ be a $(\si,D)$-vector group and let $p$ be its generic type. If $p$ is regular there is a definable subgroup
of $\gr_a$ whose generic type is regular and non-orthogonal to $p$.
\end{lem}
{\it Proof:}\\

Suppose that $G <\gr_a^d$ for some $d \in \na$. One of the projections $\pi:G \to \gr_a$ must have an infinite image
in $\gr_a$.
Let $a$ realize $p$, then $\pi(a)$ realizes the generic type of $H=\pi(G)$; this type is $tp(\pi(a)/K))$
which is also regular. Hence $H$ satisfies the conclusion of the lemma.\\
$\Box$

We have now all we need to prove \ref{teoarcs}.

{\it Proof of Theorem \ref{teoarcs}}:\\

By \ref{arc15} there is a $(\si,D)$-vector group $G$ whose generic type $q$ is regular and non-orthogonal to $p$,
 by \ref{arc16} there is a definable subgroup $H$ of the additive group  whose generic type $r$ is regular and
non-orthogonal to $q$. By transitivity of non-orthogonality on regular types, $p \not\perp r$.\\
$\Box$

\section{Partial differential fields with an automorphism}\label{sec:dncfa}

As pointed out in \cite{omar_1}, the previous work can be generalised to fields with several commuting derivations and a commuting automorphism of characteristic zero. 
We give the details.

A differential field is a field equipped with a finite set of commuting derivations
$\Delta = \{D_1, . . . , D_n\}$ . Then the ring of differential
polynomials over $K$ is defined as the ring of polynomials on the variables
$D_1^{m_1} \cdots D_n^{m_n}X$, we denote it by $K[X]_{\Delta}$
 and we can equip it with
structure of a differential ring by extending $\Delta$.

Let ${\mathcal L}_n$ be the language of differential rings with $m$ derivations, and $DF_{n}$ the
theory of differential fields of characteristic zero with $n$ commuting derivations
over ${\mathcal L}_n$.
In \cite{dcf_n}, McGrail showed that this theory has a model companion, the theory of
differentially closed fields $DCF_{n}$, and proved that it is a complete $\omega$-stable theory
which eliminates quantifiers and imaginaries. It has a noetherian topology, defined 
by zeros of ideals of differential polynomials, known as $\Delta$- topology or
Kolchin topology.

A partial differential field with an automorphism is a differential field with $n$ commuting derivations and a commuting automorphism.
In \cite{omar_1} the author showed that for the case of characteristic zero, the class of partial differential fields with an automorphsim has a model companion. We denote it by $D_nCFA$.

The axiomatisation is the following (2.1 of \cite{omar_1}):

\begin{fact}
Let $(K,\Delta, \si)$ be a differential-difference field. Then $(K,\Delta, \si)$ is existentially closed if and only if 
\begin{enumerate}
\item $(K,\Delta) $ is a model of $DCF_{n}$
\item Suppose $V$ and $W$ are irreducible $\Delta$-closed sets such that 
$W \sse V \times V^\si$
and $W$ projects $\Delta$-dominantly onto both $V$ and $V^\si$. If $O_V$ 
and $O_W$ are
nonempty $\Delta$-open sets of $V$ and $W$, respectively, 
then there is $a \in O_V$ such
that $(a, \si a) \in O_W$.
\end{enumerate}
\end{fact}

It is not known if this is a first order axiomatisation, but 2.3 of \cite{omar_1} gives us one that is first order (using characteristic sets of ideals).

As it is proved in \cite{omar_1}, the model theory of $D_nCFA$ is quite similar to the one of ${\it DCFA}$: it is not complete but its completions are easily described, they eliminate imaginaries, are not stable but are supersimple and quantifier-free stable. Section 2 of \cite{rbrank1} can be easily generalised to prove that the 
$\s$-rank of a generic of a model of $D_nCFA$ is $\omega^{n+1}$.

Now we want to prove a version of \ref{arc18} for $D_nCFA$.

Let $\mathcal U$ be a saturated model of $D_nCFA$, let $E =acl(E) \subset {\mathcal U}$, and let $G$ be a connected $\Delta$-group defined over $K=acl(K)$.
Let $G^{(n)}= G \times G^{\si} \times \cdots \times G^{{\si}^n}$ and let 
$q_n: G \to G^{(n)}$ such that $q(x)= (x, \si(x), \cdots, \si^n(x))$.

Let $g$ be a generic of $G$ such that $q_n(g)$ is $\Delta$-independent over $K$.
Then $q_n(g)$ is a ($\Delta$-)generic of $G^{(n)}$ an thus $q_n(G)$ is $\Delta$-dense in $G^{(n)}$. This 
implies that $G^{(n)}$ is connected.

Let $H$ be a definable subgroup of $G$ and let $H^{(n)}$ be the $\Delta$-Zariski closure of $q_n(H)$
in $G^{(n)}$.

Let $\tilde{ H} ^{(n)}= \{x \in G: q_n(x)  \in H^{(n)}\}$. These groups form a decreasing
sequence of quantifier-free definable subgroups of $G$ containing $H$. By Noetherianity this
sequence is finite, so there is $N$ such that
$\tilde{H}^{(N)}= \cap_n \tilde{ H} ^{(n)} $. Thus this is the  Zariski closure of $H$ in $G$.

As in \cite{rbrank1}, we can prove that  $[\tilde H:H]<\infty$.



Let $\mathcal{C}_i$ be the field
of constants of $D_i$ . The field of total constants is $\mathcal{C}= \cap \mathcal{C}_i$. Let 
$\Fix \si$ be the fixed field. Since every algebraic subgroup of a vector group is
defined by linear equations,  each
$\tilde H^n$ is defined by linear differential equations, and thus
$H$ is defined by linear $(\si,\Delta)$-equations and this implies that it is
a $(\Fix \si\cap {\mathcal C})$-vector space.
 Thus, following \cite{rbrank1}, we can prove the following lemma.

\begin{lem}\label{arccor1}\hspace{10cm}
\begin{enumerate}
\item Let $H$ be a quantifier-free definable subgroup of $\gr_a^n$. Then $H$
is a $(\Fix\si \cap {\mathcal C})$-vector space, so it is divisible and has therefore no proper subgroup of
finite index. This implies that every definable subgroup of $\gr_a^n$ is
quantifier-free definable. 
\item Let $G$ be a definable subgroup of $\gr_a^n$, and $H$ a definable subgroup of $G$. Then $G/H$ is definably isomorphic
to a subgroup of $\gr_a^l$ for some $l$.
\end{enumerate}
\end{lem}




We can apply the above lemma to definable subgroups of the additive group, and following \cite{rbrank1} again, obtain the following proposition.

\begin{prop}\label{dncfa18}
Let $G$ be a definable subgroup of $\gr_a^n$. If $G$ has infinite dimension then $\s(G)\geq \omega$.
\end{prop}

%

In \cite{omar_1}, the author proved the dichotomy for finite dimensional types.

\begin{fact}
\label{J_n}
Let $(\mathcal{U},\si, D)$ be a saturated model of $D_nCFA$ and let $K=acl(K)\subset \mathcal{U}$. Let $p\in S(K)$ be a finite dimensional type of $\s$-rank 1.
Then $p$ is either 1-based or non orthogonal to $\Fix \si \cap {\mathcal C}$
\end{fact}

Now in order to define arc spaces for a $(\si, \Delta)$-variety we must define
the differential prolongation of an algebraic variety. We use in this case the approach of \cite{arcs} (section 2).

Let $(K, \Delta=(D_1, \cdots, D_n))$ be a differential field. Define
$K_m=K[\eta_1, \cdots, \eta_n]/(\eta_1, \cdots, \eta_n)^{m+1}$.
Let $E: K \to  K_m$ defined by:

$$a \mapsto \sum_{0 \leq \alpha_1 +\cdots + \alpha_n}
\cfrac{1}{\alpha_1!  \cdots \alpha_n !}
D_1^{\alpha_1} \cdots D_n^{\alpha_n}(a)\eta_1^{\alpha_1}\cdots\eta_n^{\alpha_n}
 $$

If $V$ is an algebraic variety defined over $K$, the $m$-th $\Delta$-prolongation
$\tau_m V$ of $V$ is the Weil restriction of $V \otimes_E K_m $ from $Spec(K_m)$
to $Spec(K)$.

The reduction maps $K_l \to K _m$ for $l \geq m$ imply that the prolongations
form a projective system $\pi_{l,m}:\tau_l \to \tau_m$. If we
identify $\tau_0$ with the identity and denote $\pi_m,0$ as $\pi_m$, we
obtain the projection $\pi_m: \tau_m V \to V$.
The map $\nabla_m: V \to \tau_mV$ defined by 
$$x \mapsto \sum_{0 \leq \alpha_1 +\cdots + \alpha_n}
\cfrac{1}{\alpha_1!  \cdots \alpha_n !}
D_1^{\alpha_1} \cdots D_n^{\alpha_n}(x)\eta_1^{\alpha_1}\cdots\eta_n^{\alpha_n}
 $$

is a section of $\pi_m$ and $\nabla_m(V)$ is $\Delta$-dense in $\tau_mV$.

Using the definition of $q_m$ and $S_m(V)$ from the previous section
we define $\Phi(V)=\tau_m(S_m(V)) (= S_m(\tau_mV))$ and
$\varphi = \nabla_m \circ q_m$.

Now we can define arc spaces of $(\si, \Delta)$-varieties as before:

We extend $\si$ and $D_i$ to $K^{(m)}$ by defining $\si(\eta_j)=\eta_j$
and $D_i \eta_j=0$. Since $\si$ and $D_i$ commute, 
we can identify $\mathcal{A}_r(S_m(V))(K)$ with $S_m(\mathcal{A}_r(V))(K)$
and we can  assume that $\mathcal{A}_r(\Phi_m(V))(K) = \Phi_m(\mathcal{A}_m(V))(K)$.

If $X$ is a $(\si,\Delta)$-variety given as a $(\si,\Delta)$-closed subset of an
algebraic variety $\bar{X}$, we define $\Phi_m(X)$ as the Zariski closure
of $\psi_m(X)$ in $\Phi_m(\bar{X})$. $X$ is determined by the prolongation sequence
$\{ t_{l,m}:\Phi_l(X) \to \Phi_m(X)):l \geq m \}$. \\


We extend the previous notion of non-singular point, and as before we can show that if
 $V$ be a $(\si,D)$-variety given as a closed subvariety of an algebraic variety $\bar{V}$.
If $m \in \na$ and $a\in V(K)$ is a non-singular point.
Then $\{\mathcal{A}_m(t_{r,s}): \mathcal{A}_m\Phi_r(V)_{\psi_r(a)} \to \mathcal{A}_m\Phi_s(V)_{\psi_s(a)}, r \geq s \}$
 form the
 $(\si,D)$-prolongation sequence of
a $(\si,D)$-subvariety of $\mathcal{A}_m\bar{V}_a$.

So we define the $m$-th arc space of $V$ at $a$, $\mathcal{A}_mV_a$, to be this subvariety.

 Now we remark that in  \cite{Levin},  the author proved his theorem for a field
 with $n$ derivations and $k$ automorphisms, so for our case
 subfield $F$ of $K$ and $a \in K$ there is a numerical polynomial 
 $P_{a/F}(X) \in {\mathbb Q}[X]$ of degree at most $n+1$, such that for sufficiently large $r \in \na$,
 $P_{a/F}(r)=\td(\psi_r(a)/F)$. This implies that the notion of $(\si, \Delta)$-type
 if a $(\si, \Delta)$ -variety is well defined and it is less or equal than $n+1$.
 

If $p$ is a type over $K$, there is a finite sequence of regular types $p_1, \cdots, p_k$ such that $m(p) \geq m(p_i)$ for all $i$ and $p$
is domination-equivalent to $p_1 \times \cdots \times p_k$.
As in the previous section, we can prove that given a type $p$ over $K$ and
 a regular type $q$, such that $p \not\perp q$, there are  $c$ and $L$ such that
 the type $r=tp(c/L)$
is regular, $q \not\perp r$, $m(r) \leq m(p)$ and $\s(r)= \omega^i$ for
 $i \in \{0,1, \cdots, n+1 \}$. This is essentially what we need to prove \ref{arc13}.

The other propositions can be translated almost verbatim to the case of several derivatives, thus obtaining the following version of \ref{arc16}.

\begin{lem}\label{dncfa16}
Let $G$ be a $(\si,\Delta)$-vector group and let $p$ be its generic type. If $p$ is regular there is a definable subgroup
of $\gr_a$ whose generic type is regular and non-orthogonal to $p$.
\end{lem}
Putting \ref{dncfa18}, \ref{J_n} and \ref{dncfa16} together we obtain Zilber's dichotomy.

\begin{theorem}
\label{teo_main}
Let $({\mathcal U}, \si, \Delta)$ be a saturated model of $D_nCFA$ and ${\mathcal C}$
its field of total constants. Let $K = acl(K) \subset {\mathcal U}$ and let $p$ be a type over $K$.
If $p$ is a type of $\s$-rank 1, then it is either 1-based or non orthogonal to $\Fix \si \cap {\mathcal C}$.
\end{theorem}

 \bibliographystyle{plain}
\def\rasp{\leavevmode\raise.45ex\hbox{$\rhook$}} \def\cprime{$'$}


\begin{thebibliography}{10}

\bibitem{simple2}
S.~Buechler, A.~Pillay, and F.~Wagner.
\newblock Supersimple theories.
\newblock {\em Journal of the AMS}, 14(1):109--124, 1998.

\bibitem{rbdcfa1}
R.~Bustamante~Medina.
\newblock Differentially closed fields of characteristic zero with a generic
  automorphism.
\newblock {\em Revista de Matem\'atica: Teor\'{\i}a y Aplicaciones},
  14(1):81--100, 2007.

\bibitem{rbjets}
R.~Bustamante~Medina.
\newblock Algebraic jet spaces and {Z}ilber's dichotomy in {DCFA}.
\newblock {\em Revista de Matem\'atica: Teor\'{\i}a y Aplicaciones},
  17(1):1--12, 2010.

\bibitem{rbrank1}
R.~Bustamante~Medina.
\newblock Rank and dimension in difference-differential fields.
\newblock {\em Notre Dame Journal of Formal Logic}, 52(4):403--414, 2011.

\bibitem{salinas}
Zo{\'e} Chatzidakis and Ehud Hrushovski.
\newblock Model theory of difference fields.
\newblock {\em Trans. Amer. Math. Soc.}, 351(8):2997--3071, 1999.

\bibitem{zopi}
Zo{\'e}. Chatzidakis and Anand Pillay.
\newblock Generic structures and simple theories.
\newblock {\em Ann. Pure Appl. Logic}, 95(1-3):71--92, 1998.

\bibitem{hrusok}
E.~Hrushovski and Z.~Sokolovic.
\newblock Minimal subsets of differentially closed fields.
\newblock {\em Preprint}, 1992.

\bibitem{omar_1}
O.~Le\'on ~S\'anchez.
\newblock On the model companion of partial differential fields with an automorphism.
\newblock {\em Israel Journal of Mathematics}, 212:419–442, 2016.

\bibitem{Levin}
A.~B. Levin.
\newblock Multivariate difference-differential dimension polynomials and new
  invariants of difference-differential field extensions.
\newblock In {\em Proceedings of the 38th International Symposium on Symbolic
  and Algebraic Computation}, ISSAC '13, pages 267--274, New York, NY, USA,
  2013. ACM.

\bibitem{arcs}
R.~Moosa, A.~Pillay, and T.~Scanlon.
\newblock Differential arcs and regular types in differential fields.
\newblock {\em J. Reine angew. Math.}, 2011(620):35--54, 2008.

\bibitem{dcf_n}
T.~McGrail
\newblock The model theory of differential fields with finitely many commuting derivations.
\newblock {\em Journal of Symbolic Logic}, 65(2):885–913, 2000.


\bibitem{moosa_scanlon_2009}
R.~Moosa and T.~Scanlon.
\newblock Generalized hasse–schmidt varieties and their jet spaces.
\newblock {\em Proceedings of the London Mathematical Society. Third Series},
  2, 08 2009.

\bibitem{moosa_scanlon_2010}
R.~Moosa and T.~Scanlon.
\newblock Jet and prolongation spaces.
\newblock {\em Journal of the Institute of Mathematics of Jussieu},
  9(2):391–430, 2010.

\bibitem{jets}
Anand Pillay and Martin Ziegler.
\newblock Jet spaces of varieties over differential and difference fields.
\newblock {\em Selecta Math. (N.S.)}, 9(4):579--599, 2003.

\bibitem{wag}
Frank~O. Wagner.
\newblock {\em Simple theories}.
\newblock Kluwer Academic Publishers, Dordrecht, 2000.

\end{thebibliography}

\end{document}